\documentclass[a4paper]{article}

\title{Reconstructing directed and weighted topologies of phase-locked oscillator networks}
\date{}

\author{
Francesco Alderisio\footnotemark[2] \ , Gianfranco Fiore\footnotemark[2] \ \& Mario di Bernardo\footnotemark[2] \footnotemark[3] \ \text{*}
}

\usepackage{amsmath}
\usepackage{amssymb}
\usepackage{amsthm}
\usepackage{hyperref}
\usepackage{graphicx}
\usepackage{subfig}
\usepackage{epstopdf}
\usepackage{booktabs}
\usepackage{tabularx}  
\usepackage{xfrac}
\usepackage{color}
\usepackage{lineno}
\usepackage[margin=0.8in]{geometry}

\usepackage{doi}
\newtheorem{dfn}{Definition}

\begin{document}

\begin{titlepage}

\maketitle
The formalism of complex networks is extensively employed to describe the dynamics of interacting agents in several applications. The features of the connections among the nodes in a network are not always provided beforehand, hence the problem of appropriately inferring them often arises. Here, we present a method to reconstruct directed and weighted topologies (REDRAW) of networks of heterogeneous phase-locked nonlinear oscillators. We ultimately plan on using REDRAW to infer the interaction structure in human ensembles engaged in coordination tasks, and give insights into the overall behavior.

\footnotetext[2]{Department of Engineering Mathematics, Merchant Venturers Building, University of Bristol, Woodland Road, Clifton, Bristol BS8 1UB, United Kingdom (\texttt{enmdb@bristol.ac.uk})}
\footnotetext[3]{Department of Electrical Engineering and Information Technology, University of Naples Federico II, Via Claudio 21, 80125 Naples, Italy (\texttt{mario.dibernardo@unina.it})}

\end{titlepage}

\setcounter{page}{2}

\section{Introduction}
The study of complex networks has been a subject of great interest \cite{strogatz2001exploring,albert2002statistical,ravasz2003hierarchical,boccaletti2006complex}. Indeed, networks provide a rigorous formalism to describe phenomena involving populations of interacting agents in fields as diverse as Physics, Engineering, Biology, Chemistry, Social Science and the Internet \cite{buchanan2003nexus,barabasi2004network,arenas2008synchronization,cardillo2014evolutionary}.

Even though the topology of the interactions among the nodes is often provided, there might be cases in which there is no \emph{a-priori} knowledge of the network structure. Instead, the topology of the interconnections among the nodes has to be reconstructed (or reverse-engineered) from data. Specifically, given a set of simultaneously recorded time series, it is important to identify the information structure characterized by hidden dependencies between the components of a complex dynamical system \cite{lungarella2007information}. In particular, it may be required to identify causal dependencies among the nodes or quantify the flow of information across a specific network. Such flow can be \emph{unidirectional} (the behavior of node/agent A affects that of node/agent B, but not vice versa) or \emph{bidirectional} (the behavior of each of the two nodes/agents depends on that of the other).
Moreover, the strength by which the agents influence each other might be different; the information flow in this case is therefore said to be \emph{weighted}.
Knowing the topology of the interactions among the agents in a network would allow for a better understanding of their behavior in order to analyze it and predict it, design networks with pre-described functions, and potentially control them \cite{liu2011controllability,menolascina2014vivo}.

The problem of inferring the existence and direction of links among nodes in a network from a given data-set is a reverse engineering task often referred to as \emph{network reconstruction} or \emph{network inference}. It is of crucial importance to determine which of the inferred connections correspond to links existing in the real interaction structure (\emph{true positives}), and which instead represent only functional connectivities between nodes (\emph{false positives}) \cite{kralemann2011reconstructing}. This problem arises in several contexts, e.g., reconstructing functional activation/repression links in gene regulatory networks \cite{crampin2004mathematical,bansal2007infer,barzel2013network} or causal relationships in stochastic processes \cite{sun2015causal}, understanding the structure of social interactions in a group from communication data \cite{bagrow2005local,sun2014causation,villaverde2014reverse}, and inferring functional relationships between areas in the brain from EEG data \cite{staniek2008symbolic} or in physiological systems \cite{schreiber2000measuring}. However, reviews of network inference methods have found large discrepancies among the results of different algorithms \cite{de2010advantages,marbach2010revealing}.

A notable case is that of networks of nodes whose dynamics is oscillatory (e.g., neurones, cellular cycles, synthetic biological oscillators, groups of walking autonomous robots). For these networks, methods for reconstructing their structure from data on their dynamics have been proposed. Examples include the work in \cite{yu2006estimating}, dealing with the problem of inferring directed and weighted topologies, as well as that in  \cite{arenas2006synchronization}, investigating community detection in undirected and unweighted networks of Kuramoto oscillators \cite{kuramoto2012chemical}, and those in \cite{kralemann2011reconstructing,timme2007revealing}, tackling the reconstruction of directed and unweighted topologies of small networks of coupled limit cycle oscillators. In these works, none of the inferred connections are removed, or if they are, no formal method is presented on how to choose the thresholds according to which false positive links are to be cut-off.

Here, we present REDRAW (REconstruction of DiRected And Weighted topologies) as a significant extension of the methodology firstly proposed in \cite{arenas2006synchronization}, originally conceived to only detect communities in undirected and unweighted networks of Kuramoto oscillators, widely used to describe synchronization phenomena within populations of interacting agents \cite{acebron2005kuramoto,antonioni2016coevolution}. We show that REDRAW allows to effectively infer directed and weighted links among a group of coupled oscillators, and present an algorithm that allows to set the cut-off thresholds for removing false positives. We validate the method on a number of representative examples, including a set of real-world networks obtained from \cite{villaverde2014mider,basso2005reverse,diebold2014network,vinayagam2011directed,chang2011performance}.

\section{Materials and Methods}

The mathematical model used to validate REDRAW is a network of $n$ nonuniform Kuramoto oscillators \cite{dorfler2012synchronization,zhang2016bounded}, described by

\begin{equation}
\label{eqn:1}
\dot{\theta_i} = \omega_i + \frac{c}{n} \sum_{j=1}^n a_{ij} \sin(\theta_j-\theta_i-\phi_{ij}), \ i=1,2,\ldots,n
\end{equation}
where $\theta_i \in [-\pi,\pi]$ represents the phase of the $i$th oscillator, $\omega_i>0$ its natural frequency, $c > 0$ the global coupling strength among all nodes in the network, and $a_{ij}$ the local influence that node $j$ has on node $i$. For directed and weighted topologies, in general $a_{ij} \neq a_{ji}$ with $a_{ij} \ge 0 \ \forall i,j$. Moreover, the phase shift

\begin{equation}
\label{eqn:2}
\phi_{ij}:=\begin{cases}
\frac{\phi}{a_{ij}}, & a{ij}>0 \\
0, & a_{ij} = 0
\end{cases}
\end{equation}
represents how much node $i$ lags behind node $j$, with $\phi \in [0,\frac{\pi}{2}]$. Note that the higher the influence $a_{ij}$ that node $j$ has on node $i$, the lower the phase shift $\phi_{ij}$ (i.e., the less node $i$ lags behind node $j$).

We start from a data-set consisting of $n$ time series of duration $T>0$, one for each of the $n$ oscillators in the network of interest, and assume that $K$ experiments are available for a total of $K \times n$ time series of length $T$. The proposed method aims at inferring the topology of the interactions among the nodes by estimating the influence that each agent in the network has on the others, provided \emph{phase-locking} is achieved by the oscillators \cite{suppinfoPL}.

REDRAW is a method consisting in a sequence of six steps: first the relative phases among all the nodes are evaluated and parameters representing their respective correlation are derived from them (\emph{Step 1}--\emph{Step 4}), then two filtering steps (\emph{Step 5} and \emph{Step 6}) are carried out to remove possible false positives in the interactions among the nodes. More specifically:

\begin{enumerate}
\item[\emph{Step 1}.] The relative phase $\Delta\theta_{ij,k}(t)=\theta_{i,k}(t)-\theta_{j,k}(t)$, $\forall k=1,\ldots,K$ with $t \in [0,T]$ is evaluated, where negative values of $\Delta\theta_{ij,k}(t)$ indicate that node $i$ lags behind node $j$ at time instant $t$ of the $k$th experiment.
\item[\emph{Step 2}.] The parameter $\zeta_{ij,k}(t)$ defined as:
\begin{equation}
\label{eqn:3}
\zeta_{ij,k}(t):=\begin{cases}
\frac{1+\cos[\Delta\theta_{ij,k}(t)]}{2}, & \Delta\theta_{ij,k}(t)\le 0 \\
0, & \Delta\theta_{ij,k}(t)> 0
\end{cases}
\end{equation}
is calculated for all the possible pairs of nodes in the network. Note that $\zeta_{ij,k}(t) \neq \zeta_{ji,k}(t)$ and that $\zeta_{ij,k}(t) \in [0,1]$, where $\zeta_{ij,k}(t)=1$ represents the maximum level of influence that node $j$ has on node $i$. On the other hand, when $\Delta\theta_{ij,k}(t)>0$, agent $j$ lags behind agent $i$, hence the influence that the former has on the latter is assumed to be negligible and set to 0.
\item[\emph{Step 3}.] The average value over time of the parameter defined in Eq. \eqref{eqn:3} is evaluated as:
\begin{equation}
\label{eqn:4}
\rho_{ij,k} := \frac{1}{T} \int_0^T \zeta_{ij,k}(t) \ dt
\end{equation}
where $T$ is the duration of each time series.
\item[\emph{Step 4}.] The link representative of the overall influence that node $j$ has on node $i$ is inferred by averaging the parameter defined in Eq. \eqref{eqn:4} over the total number $K$ of experiments:
\begin{equation}
\label{eqn:5}
\rho_{ij} := \frac{1}{K} \sum_{k=1}^K \rho_{ij,k}
\end{equation}

\item[\emph{Step 5}.] Data Processing Inequality (DPI): the interactions $\rho_{ij}$ among all triplets of connected nodes are checked, and on the basis of their intensity one of them is possibly regarded as a false connection \cite{margolin2006aracne}. Suppose to have, for instance, three connected nodes $w$, $y$ and $z$: if the influence $\rho_{zw}$ that $w$ has on $z$ is lower than both the one that $w$ has on $y$, that is $\rho_{yw}$, and the one that $y$ has on $z$, that is $\rho_{zy}$, then the link between the pair $(w,z)$ is removed, that is $\rho_{zw}$ is set to $0$. In particular, in order not to remove links in triplets that could actually be connected in the original topology, a further condition is added to the standard DPI, according to which false positives $\rho_{ij}$ are removed as long as their value is below a certain threshold $\nu$, with $0\le \nu<1$ (Fig. \ref{fig:1}). The higher the value of $\nu$, the less connected triplets are found in the network.

\begin{figure}[h!]
\centering
\includegraphics[scale=1.25]{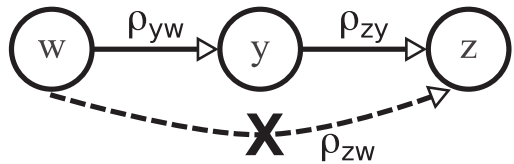}
\caption{\label{fig:1} Data Processing Inequality. An edge going out of node $i$ and coming in node $j$, where the direction is given by the arrow, is representative of the fact that node $j$ is influenced by node $i$. Parameter $\rho_{zw}$ is set to $0$ as long as the three following conditions are verified simultaneously: 1) $\rho_{zw} < \rho_{yw}$, 2) $\rho_{zw} < \rho_{zy}$, 3) $\rho_{zw} < \nu$ , where $\nu$ is a threshold value.}
\end{figure}

\item[\emph{Step 6}.] Network Thresholding: all parameters $\rho_{ij}$ whose value is below a certain threshold $\mu$ are set to $0$, with $0\le \mu \le \nu$ (the higher the value of $\mu$, the sparser is the reconstructed network structure).
\end{enumerate}

In order to investigate how the topology of the interactions among the nodes evolves over time, the time interval $[0,T]$ can be partitioned in $L$ time windows $[t_l,t_{l+1}]$ of length $\Delta T_l$, with $l=0,\ldots,L-1$.
REDRAW can thus be applied to measurements available for each of the $L$ time intervals, so that the corresponding network structure can be associated with each of them.
The choice of $\Delta T_l$ is dependent on the specific application of interest and can be selected heuristically by trial-and-error.

Standard metrics are used to assess the performance of REDRAW \cite{dasgupta2016models}: Positive Predicted Value (PPV), Accuracy Rate (ACC), True Positive Rate (TPR) and False Positive Rate (FPR). They are provided as percentages, so that a network topology is perfectly inferred when $PPV=ACC=TPR=100 \%$ and $FPR=0 \%$. Details on their definition are given in \cite{suppinfo}.

Notably, we provided an algorithm for a sensible choice of the values that the thresholds $\nu$ and $\mu$ should take to remove possible false positives from the interconnections inferred in \emph{Step 4}.
Specifically, such algorithm relates topologies of different size to corresponding ranges of thresholds that maximize the metrics here considered when reconstructing an unknown network of interest. In particular, this is an a-priori method whose input consists in the size $n$ of the network to reconstruct and $N$ random graphs (each made up of $n$ nodes) following the Erd{\"o}s-R{\'e}nyi model \cite{erdos1959random,erdos1960evolution}, and whose output is an acceptable range of values for $\nu$ and $\mu$.
In particular, the algorithm reconstructs the assigned random graphs for different values of $\nu$ and $\mu$, each leading to different values of the performance metrics, respectively.
Only pairs ($\nu$,$\mu$) giving rise to acceptable values for such metrics are provided as output and should then be selected to reconstruct unknown topologies of interest.
For further details on the aforementioned algorithm see \cite{suppinfo}.

\section{Results}
\begin{figure}[h!]
\centering
\includegraphics[scale=.8]{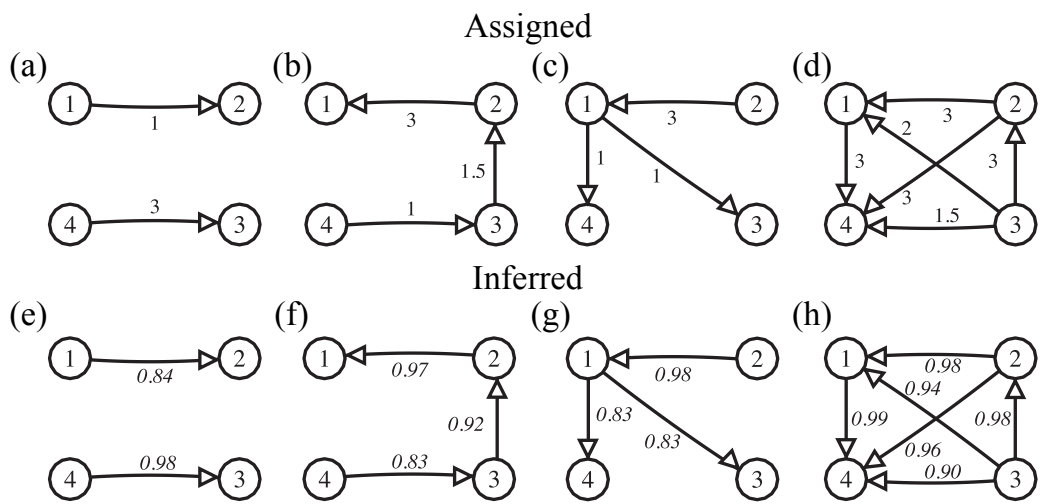}
\caption{\label{fig:2} Assigned and inferred topologies, $n =  4$. The topologies on top (a-d) represent those used in the numerical simulations to generate the data-set then employed to obtain the inferred topologies, respectively represented in the bottom panels (e-h). The numerical values on the edges of the assigned topologies represent the values of $a_{ij}$ in the model described in Eq. \eqref{eqn:1}, whereas the italic numerical values on the edges of the inferred topologies represent parameters $\rho_{ij}$ estimated by REDRAW.}
\end{figure}

A network of $n =  4$ nonuniform Kuramoto oscillators described by Eq. \eqref{eqn:1} was first considered. Four different topologies (assumed to be unknown), respectively represented in Figs. \ref{fig:2}(a)-\ref{fig:2}(d), were investigated. For each topology, $K = 50$ experiments of duration $T = 30$s were numerically simulated. For each node, experiment and topology, the natural frequencies $\omega_i$ were randomly extracted from the interval $[1,2]$ rad s$^{-1}$, and so were the initial conditions $\theta_i(0)$ from the interval $[-\pi,\pi]$, $i=1,\ldots,4$. The model parameters were set to $c=10$ and $\phi=\frac{\pi}{4}$ so that phase-locking could be achieved, while the thresholds were set to $\nu=0.9$ and $\mu=0.8$ following the algorithm in \cite{suppinfo}.

The reconstructed networks are represented in Figs. \ref{fig:2}(e)-\ref{fig:2}(g), respectively. For each topology, neither missing links nor false positives are found, and the directionality of the links is correctly inferred for all of them ($PPV=ACC=TPR=100 \%$, $FPR=0 \%$). As for the weights, their magnitudes' relative relationship in the assigned topologies is correctly inferred as well. For instance, note how $a_{34} < a_{23} < a_{12}$ in Fig. \ref{fig:2}(b) corresponds to $\rho_{34} < \rho_{23} < \rho_{12}$ in Fig. \ref{fig:2}(f), or how $a_{31} = a_{41} < a_{12}$ in Fig. \ref{fig:2}(c) corresponds to $\rho_{31} = \rho_{41} < \rho_{12}$ in Fig. \ref{fig:2}(g). 

\begin{figure}[ht]
\centering
\includegraphics[scale=0.8]{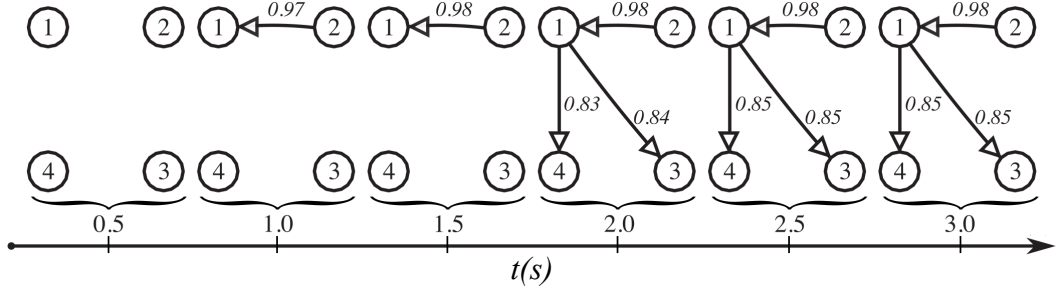}
\caption{\label{fig:3} Topology of Fig. \ref{fig:2}(c) inferred every $0.5$s. The italic numerical values on the edges represent parameters $\rho_{ij}$ estimated by REDRAW.}
\end{figure}

The evolution over time of the reconstructed network with topology represented in Fig. \ref{fig:2}(c) was inferred over time windows of length $\Delta T_l = \Delta T = 0.5$s, and it is shown in Fig. \ref{fig:3} (analogous results are found for the other topologies, data not shown). Notably, none of the nodes are connected before $t = 0.5$s, and then the stronger connection between nodes $1$ and $2$ is inferred before the others (the reconstructed topology does not change after $t = 3$s).

\begin{figure}[ht]
\centering
\includegraphics[scale=0.8]{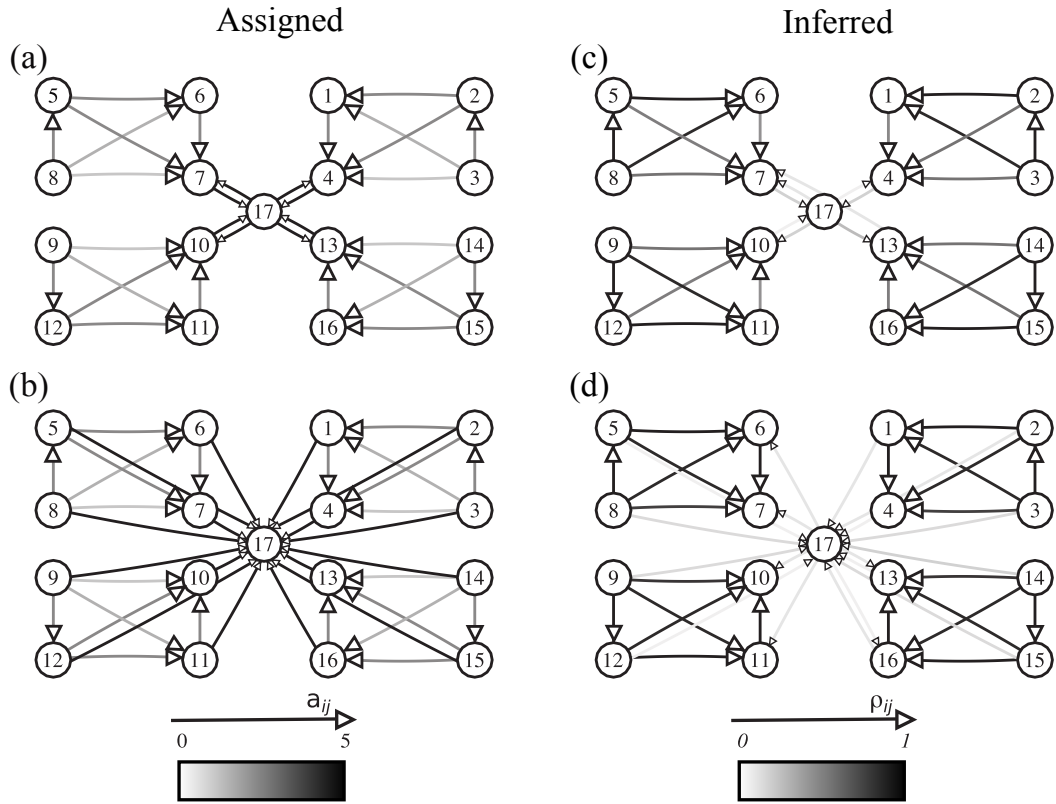}
\caption{\label{fig:4} Assigned and inferred topologies, $n =  17$. The geometric graph (a) and the Ravasz-Barab{\'a}si network (b) on the left-hand side represent the topologies used in the numerical simulations to generate the data-set then employed to obtain the inferred ones, respectively depicted on the right-hand side (c,d). Different scales of gray quantify the numerical value of $a_{ij}$ for the assigned topologies, and $\rho_{ij}$ estimated by REDRAW for the inferred ones.}
\end{figure}

Next, we tested REDRAW on a larger network of $n=17$ nodes. We generated data numerically by simulating the network over two different topologies (assumed to be unknown) with $c = 40$ so that phase-locking could be achieved, and set the threshold $\mu=0.35$ while keeping all other parameters as in the previous example.
The structures being considered were obtained by interconnecting four sub-networks with the structure represented in Fig. \ref{fig:2}(d), either through a central hub in a  \emph{geometric graph} configuration \cite{dall2002random} shown in Fig. \ref{fig:4}(a), or as the \emph{Ravasz-Barab{\'a}si graph} \cite{ravasz2003hierarchical} shown in Fig. \ref{fig:4}(b), respectively.

The directionality of the geometric graph represented in Fig. \ref{fig:4}(a) is correctly inferred for all its edges [Fig. \ref{fig:4}(c)], with the only exception of a missing link from node $13$ to node $17$, replaced in the reconstructed topology by a link from node $13$ to node $7$. Although in the assigned topology the values $a_{ij}$ of the edges within each sub-topology are lower than those connecting them to the central hub, the opposite result is found for $\rho_{ij}$ in the inferred topology. This is due to the fact that the oscillators reach synchronization within each corresponding sub-group of four nodes (given the higher number of connections), hence exhibit greater phase mismatch with respect to the central hub.

Similar observations can be made for the Ravasz-Barab{\'a}si network represented in Fig. \ref{fig:4}(b). The directionality is correctly inferred for all the edges belonging to each of the four sub-topologies, and for $12$ out of the $16$ edges connecting them to the central hub [Fig. \ref{fig:4}(d)].
However, despite the values of $a_{ij}$ in each sub-topology being lower than those related to the links connected to the central hub, the opposite result is found for $\rho_{ij}$ in the inferred topology. This is due to the fact that, according to the assigned topology, the central hub is influenced by all the nodes in the network, hence its phase mismatch is minimized only with respect to an average value of the phases of all the other nodes.

\begin{figure}[ht]
\centering
\includegraphics[scale=0.8]{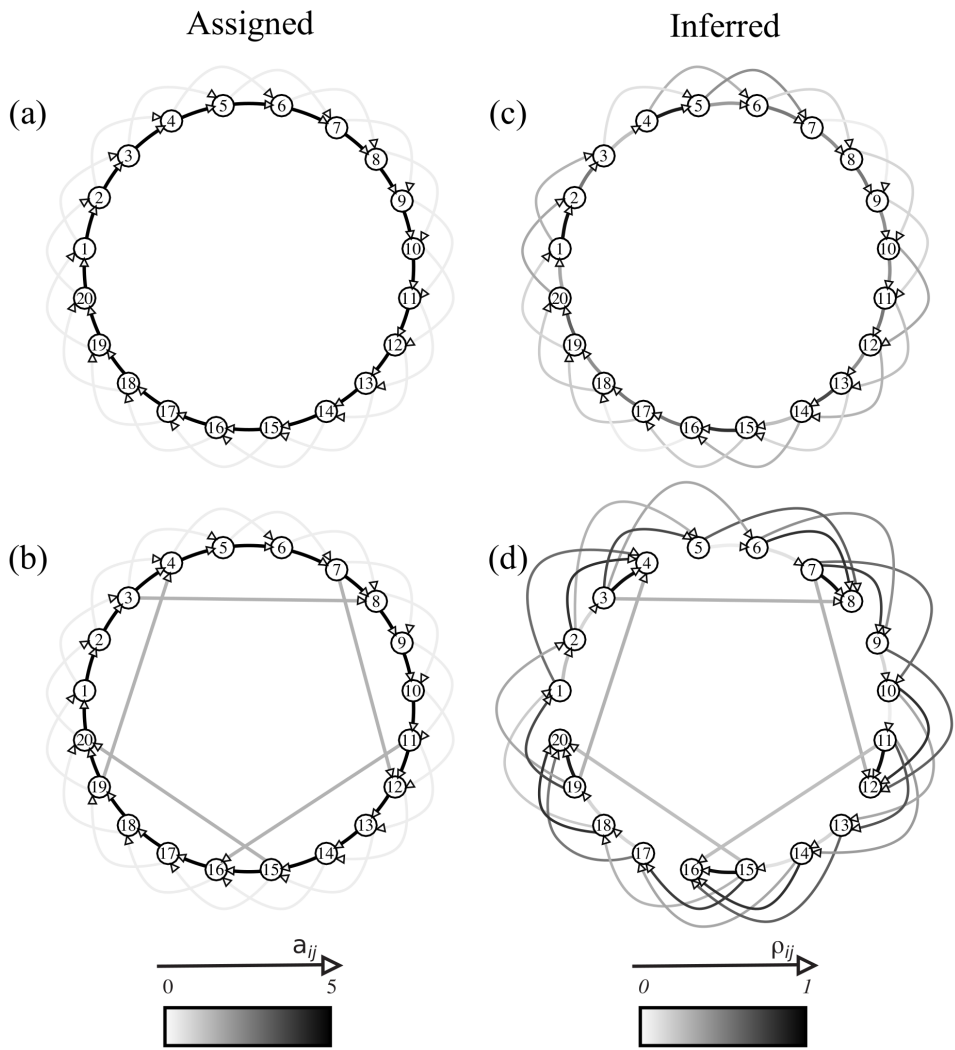}
\caption{\label{fig:5} Assigned and inferred topologies, $n =  20$. The regular network (a) and that obtained by some long-distance rewiring (b) on the left-hand side represent the topologies used in the numerical simulations to generate the data-set then employed to obtain the inferred ones, respectively represented on the right-hand side (c,d). Different scales of gray quantify the numerical value of $a_{ij}$ for the assigned topologies, and $\rho_{ij}$ estimated by REDRAW for the inferred ones.}
\end{figure}

We then applied REDRAW to reconstruct a network of $n=20$ nodes. Data were obtained via numerical simulations carried out by considering the regular network structure shown in Fig. \ref{fig:5}(a), and the structure obtained by some long-distance rewiring shown in Fig. \ref{fig:5}(b), both assumed to be unknown. The coupling strength among nodes was set as $c = 50$ in the simulations so that phase locking could be achieved. REDRAW was parameterized by setting the thresholds to $\nu=0.65$ and $\mu=0.60$, following the algorithm in \cite{suppinfo}, with the other parameters being selected as before.

The directionality of the regular network represented in Fig. \ref{fig:5}(a) is correctly inferred for all its edges [Fig. \ref{fig:5}(c)], with the only exception of a missing link from node $2$ to node $4$. For each $i$th node, $i=1,2,\ldots,20$, the edge coming from agent $i-1$ has a higher inferred weight than that coming from agent $i-2$ (note that agents $-1,0$ correspond to agents $19,20$, respectively), thus reproducing well the interactions assigned in the original topology.

The introduction of $5$ edges [Fig. \ref{fig:5}(b)] is well captured by the topology inferred in Fig. \ref{fig:5}(d), leading to the formation of an equal number of clusters [nodes $(1-4)$, $(5-8)$, $(9-12)$, $(13-16)$ and $(17-20)$]. For each of the $5$ nodes being influenced by one of the additional links (nodes $4, 8, 12, 16, 20$), no outgoing links are inferred. This is a result of the model described in Eq. \eqref{eqn:1}: for instance, the phase of node $4$ is lowered by the influence of node $19$, thus leading to higher mismatches of the former with nodes $5$ and $6$ (similar reasoning can be carried out for nodes $8$, $12$, $16$ and $20$).  

Quantitative details on the effectiveness of the proposed method are given in Table \ref{table:1} for the topologies represented in Figs. \ref{fig:4} and \ref{fig:5}, respectively. In all the cases but for the network in Fig. \ref{fig:5}(d), PPV, ACC and TPR take values close to $100 \%$, while FPR takes values close or equal to $0 \%$. In the remaining case, despite PPV and TPR taking values around $67 \%$, the method still provides a low value of FPR ($4.5 \%$) and a high value of ACC ($92.1 \%$).

 \begin{table}
 \centering
 \caption{\label{table:1} Metrics used to validate the effectiveness of REDRAW on networks of $n=17$ and $n=20$ nodes.}
 \begin{tabular}{ccccc}
 \hline
 \textbf{Topologies} & $\mathbf{PPV (\%)}$ & $\mathbf{ACC (\%)}$ & $\mathbf{TPR (\%)}$ & $\mathbf{FPR (\%)}$ \\
 \hline \hline
Geometric graph & 96.9 & 99.3 & 96.9 & 0.04 \\
Ravasz-Barab{\'a}si & 86.4 & 97.1 & 95.0 & 2.60 \\
Regular network & 100 & 99.7 & 97.5 & 0 \\
Rewired network & 66.7 & 92.1 & 66.7 & 4.50 \\
\hline
 \end{tabular}
 \end{table}

Finally, we used REDRAW to reconstruct some real-world networks existing in the current literature \cite{villaverde2014mider,basso2005reverse,diebold2014network,vinayagam2011directed,chang2011performance}.
The values of the four metrics are detailed in Table \ref{table:2} for different topologies characterized by a different number of nodes $n$ and edges $e$, together with the value of the coupling strength $c$ in the model described in Eq. \eqref{eqn:1}, selected so that phase-locking could be achieved \cite{suppinfoPL}, and that of the thresholds $\nu$ and $\mu$, selected as suggested by the proposed algorithm \cite{suppinfo} -- note how such thresholds converge towards each other as the number of nodes increases.
All other parameters were kept as in the previous examples.
The results in Table \ref{table:2} confirm the effectiveness of REDRAW as method for reconstructing nonlinear heterogeneous oscillator networks. Indeed, for any of the considered real-world topologies, at least three out of four metrics take values representative of a correct inference. 

\begin{table}[t]
\centering
 \caption{\label{table:2} Validation of REDRAW on the reconstruction of real-world networks. This table shows the four performance metrics (in percentages) obtained for each topology, together with their own number of nodes $n$ and edges $e$, respectively. The values of coupling strength $c$ and thresholds $\nu$ and $\mu$ employed in the numerical simulations are detailed as well.}
\center
 \begin{tabular}{cccccccccc}
 \hline
 \textbf{Topologies} & $\mathbf{n}$ & $\mathbf{e}$ & $\mathbf{c}$ & $\mathbf{\nu}$ & $\mathbf{\mu}$ & $\mathbf{PPV}$ & $\mathbf{ACC}$ & $\mathbf{TPR}$ & $\mathbf{FPR}$ \\
 \hline \hline
Enzyme-Catalyzed reaction pathway \cite{villaverde2014mider} & 8 & 7 & 20 & 0.90 & 0.75 & 100 & 100 & 100 & 0 \\
Songbird brain \cite{basso2005reverse} & 12 & 13 & 35 & 0.90 & 0.65 & 63 & 93 & 77 & 5 \\
Bank stocks connections \cite{diebold2014network} & 16 & 32 & 40 & 0.80 & 0.70 & 56 & 87 & 16 & 2 \\
Human PPI \cite{vinayagam2011directed} & 23 & 22 & 40 & 0.80 & 0.50 & 100 & 98 & 45 & 0 \\
Human PPI \cite{vinayagam2011directed} & 25 & 27 & 50 & 0.80 & 0.65 & 56 & 96 & 44 & 2 \\
Hainan Power Grid Company \cite{chang2011performance} & 48 & 63 & 150 & 0.60 & 0.60 & 52 & 97 & 62 & 2 \\
\hline
\end{tabular}
 \end{table}
 
For all the topologies considered in this work, we quantified the performance metrics for different numbers of experiments and observed how, for a sufficiently high number of repetitions, REDRAW is not sensitive to the specific value of $K$. For the sake of brevity, we here reported results only for the case $K=50$; in fact, when $K\ge50$, possible fluctuations of the metrics are negligible. For more details on the performance of REDRAW for different number of experiments see \cite{suppinfoTrials}.

\section{Conclusion and Discussion}
We presented REDRAW as a network reconstruction method to infer directed and weighted links among groups of heterogeneous phase-locked nonlinear oscillators, whose interconnections were assumed to be unknown, and tested it in the particular case of networks of nonuniform Kuramoto oscillators for representative as well as real-world topologies. We observed how both directionality and weights could be correctly estimated with acceptable values of the performance metrics even in the case of larger networks.

Notably, we provided an algorithm for a sensible choice of the thresholds $\nu$ and $\mu$ employed to remove possible false positives from the inferred interconnections. Such algorithm is based on the reconstruction of random graphs following the Erd{\"o}s-R{\'e}nyi model \cite{erdos1959random,erdos1960evolution}, and is presented in details in the accompanying supplementary material \cite{suppinfo}.

We plan on using REDRAW to infer the interactions underlying small-scale human ensembles engaged in multiplayer coordination tasks that require the group members to generate sinusoidal-like motions \cite{alderisio2016entrainment,alderisio2016interaction,alderisio2016study}, and therefore give insights into their behavior when interacting together.
Specifically, from the directionality of the inferred links it is possible to detect the emerging leader in the group as the agent that influences the others the most in terms of highest number of outgoing edges \cite{d2012leadership}, and hence decide what agent(s) could be possibly entrained (e.g., by an external audio/visual signal) for steering the whole group towards a desired behavior \cite{wang2013exponential}.
Moreover, from the weights of the inferred links it is possible to gain information on the partner(s) each participant prefers to interact with and on those s/he tends to avoid. This can offer a criterion to redefine the topology of the group interactions in order for each participant to maximize synchronization with all the others, which is fundamental when a high level of coordination is required in a human ensemble, as in music \cite{badino2014sensorimotor,volpe2016measuring} and sports \cite{wing1995coordination,silva2016practice}.

\section*{Acknowledgments}
The authors wish to acknowledge support from the European Project AlterEgo FP7 ICT 2.9 -- Cognitive Sciences and Robotics, Grant Number 600610.


\clearpage

\section*{SUPPLEMENTARY MATERIAL}

\subsection*{Phase-locking}

In order for REDRAW to correctly infer the structure of interconnections among $n$ coupled oscillators, it is necessary for them to have achieved \emph{phase-locking} as defined below.

\begin{dfn}
\label{dfn:pl}
Denoting with $\theta_i(t), \ i=1,...,n$, the phase of the $i$-th oscillator in the network at time $t$, and with $r$ and $\psi$ the Kuramoto order parameters defined as:

\begin{equation}
r(t) e^{j\psi(t)} := \frac{1}{n} \sum_{i=1}^{n} e^{j\theta_i(t)} \ ,
\end{equation}
we say that \emph{phase-locking} is achieved when

\begin{equation}
c_v := \frac{\sigma}{\eta} \le \chi \ ,
\end{equation}
where $c_v$, $\sigma$ and $\eta$ are coefficient of variation, standard deviation and mean over time of $\psi(t) \ \forall t\ge \hat{t}$, respectively, and $\chi>0$ represents a certain upper bound.
\end{dfn}

As $\psi$ represents the average angular velocity of all the oscillators, the previous condition guarantees that all their phase differences are bounded $\forall \ t \ge \hat{t}$. 

\subsection*{Algorithm to select thresholds $\nu$ and $\mu$}
\begin{figure}[!h]
 \centering
\includegraphics[scale=.75]{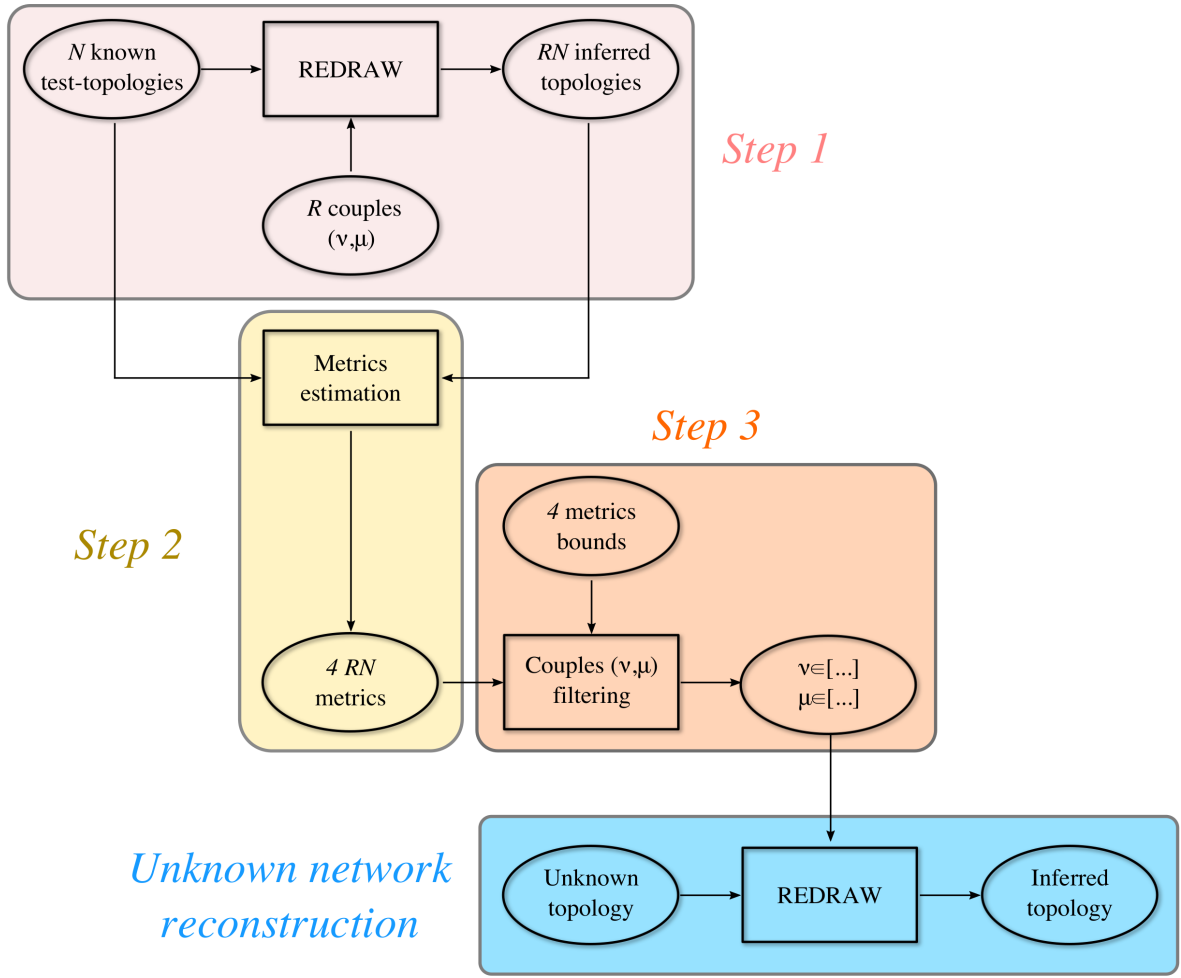}
 \caption{Algorithm to select acceptable threshold values for $\nu$ and $\mu$. \emph{Step 1}: $N$ known test-topologies are reconstructed by REDRAW for $R$ possible combinations of threshold couples $(\nu,\mu)$. \emph{Step 2}: for each of the $RN$ reconstructed topologies, four standard metrics are evaluated. \emph{Step 3}: averages of the $RN$ values for each of the four metrics are evaluated across the $N$ different structures, and thresholds $\nu$ and $\mu$ are selected so that four bound conditions are simultaneously verified. The thresholds thus obtained can then be selected to reconstruct an unknown topology of interest.}
 \label{fig:S1}
\end{figure}

An important step in REDRAW when reconstructing the interactions among the nodes of a network whose topology is unknown, is the use of two filtering thresholds $0\le \nu<1$ and $0\le \mu \le \nu$: the former is used when applying \emph{Data Processing Inequality} (DPI), the latter when removing some of the inferred links. The choice of such thresholds can be crucial, hence it is necessary to provide an \emph{a-priori} criterion according to which acceptable values can be selected for them.
Here we detail a possible algorithm to select the values of $\nu$ and $\mu$ (Figure \ref{fig:S1}).

Specifically, we suppose that a data-set is available for $N$ different network configurations (assumed to be known) of $n$ nodes reaching phase-locking, where $n$ also corresponds to the size of the unknown topology of interest. For each of the $N$ configurations, we assume $K$ experiments are available, each of duration $T$. For each experiment, time-series for all the nodes are available.
Then, we propose the following steps to select thresholds $\nu$ and $\mu$ that should eventually be employed when inferring the structure of interactions among the agents for which only experimental data is available (no a-priori information is given on their topology).

\begin{itemize}
\item[\emph{Step 1}.] Consider a two-dimensional grid made up of $R$ points individuated by $R$ respective pairs $(\nu,\mu)$. Each of the $N$ test-topologies is reconstructed by REDRAW for all the $R$ threshold couples $(\nu,\mu)$ such that $0 \le \mu \le \nu < 1$. As a result, $RN$ topologies are inferred.
\item[\emph{Step 2}.] For each of the $RN$ reconstructed topologies, four standard metrics are computed. 
Specifically, the parameter $\rho_{ij}$ inferred $\forall i,j$ through REDRAW is said to be a: True Positive (TP), if $\rho_{ij}>0$ and $a_{ij}>0$; False Positive (FP), if $\rho_{ij}>0$ and $a_{ij}=0$; True Negative (TN), if $\rho_{ij}=0$ and $a_{ij}=0$; False Negative (FN), if $\rho_{ij}=0$ and $a_{ij}>0$, 
where $a_{ij}$ is the corresponding value in the topology originally assigned.
Denoting with $N_{TP}$, $N_{FP}$, $N_{TN}$, $N_{FN}$ and $N_{TOT}:=n(n-1)$ the total number of true positives, false positives, true negatives, false negatives and possible links among all the nodes in the network, respectively, the performance metrics employed here are:
\begin{equation*}
PPV:=\frac{N_{TP}}{N_{TP}+N_{FP}}, \ ACC:=\frac{N_{TP}+N_{TN}}{N_{TOT}}, \
TPR:=\frac{N_{TP}}{N_{TP}+N_{FN}}, \ FPR:=\frac{N_{FP}}{N_{FP}+N_{TN}}.
\end{equation*}

\item[\emph{Step 3}.] For each metric, acceptable bounds are defined according to the level of accuracy of the reconstruction that is desired. We term these bounds $PPV^*$, $ACC^*$, $TPR^*$ and $FPR^*$. Averages of the $RN$ values for each of the four metrics are computed across the $N$ different structures. Thresholds $\nu$ and $\mu$ (i.e., the final output of the algorithm) are then selected so that the following conditions are simultaneously verified:
\begin{equation}
\label{eqn:metricsBounds}
E[PPV] \ge PPV^{*}, \quad E[ACC] \ge ACC^{*}, \quad E[TPR] \ge TPR^{*}, \quad E[FPR] \le FPR^{*}
\end{equation}
where $E[M]$ denotes the average value of metric $M$.
\end{itemize}

In our work we take the following choices:

\begin{itemize}
%
\item We set $\chi=35 \%$, $\hat{t}=20$s, $T=30$s, $R=5050$ (corresponding to sampling $\nu$ and $\mu$ in the interval [0,0.99] each with a step-size of $0.01$), $N=100$, $K=10$, $a_{ij}=1,0$ according to whether node $i$ is influenced by node $j$ (there exists a link going from node $j$ to node $i$) or not, respectively, $\phi=\frac{\pi}{4}$ and $c=2.5 n$ so that phase-locking could be achieved. For each node, experiment and test-topology, the initial conditions $\theta_i(0)$ are randomly extracted from the interval $[-\pi,\pi]$, and so are the oscillation frequencies $\omega_i$ from the interval $[1,2]$ rad s$^{-1}$.
\item We employ directed random graphs as test-topologies. Specifically, the Erd{\"o}s-R{\'e}nyi $G(n,p)$ model is used, where the probability $p$ of a link connecting any two nodes in the network is independent on that of the others. We set such probability to $p=\frac{\ln{n}}{2n}$, which is likely to provide weakly connected graphs (their undirected version is connected).
\item As the main goal is that of providing acceptable ranges within which $\nu$ and $\mu$ should take values for inferring unknown networks rather than known test-topologies, the metrics bounds for reconstruction of the latter need not be strict and can be chosen to be flexible. In particular, we set such bounds to PPV$^{*}=$TPR$^{*}=40 \%$, ACC$^{*}=70 \%$ and FPR$^{*}=30 \%$.
\end{itemize}

In what follows we illustrate the algorithm by using a collection of representative data-sets obtained by simulating $N=100$ directed random graphs. We repeat the application of the algorithm to networks of different sizes to illustrate the effects of varying the number of nodes $n$ in the network.

\subsubsection*{Test-topologies in \emph{Step 1}}
The second eigenvalue $\lambda_2(L)$ of the Laplacian matrix of the undirected version of the test-topologies, together with the respective coefficient of variation $c_v$ obtained in the simulations, are detailed in Table \ref{table:S1}.

\begin{table}[h]
 \caption{\label{table:S1} Values of $\lambda_2(L)$ and $c_v$ averaged over the $N=100$ directed random graphs employed as known test-topologies. For both measures, their value is represented for increasing number of nodes $n$, respectively.}
\center
 \begin{tabular}{ccccccc}
 \hline
 $n$ & 5 & 10 & 15 & 20 & 25 & 30 \\
 \hline \hline
$\lambda_2(L)$ & 0.86 & 0.46 & 0.43 & 037 & 0.40 & 0.51 \\
$c_v (\%)$ & 1.3 & 4.3 & 7.6 & 8.4 & 11 & 10 \\
\hline
\end{tabular}
 \end{table}

\clearpage
\subsubsection*{\emph{Step 2}: metrics evaluation}

\begin{figure}[!h]
 \centering
\includegraphics[scale=.9]{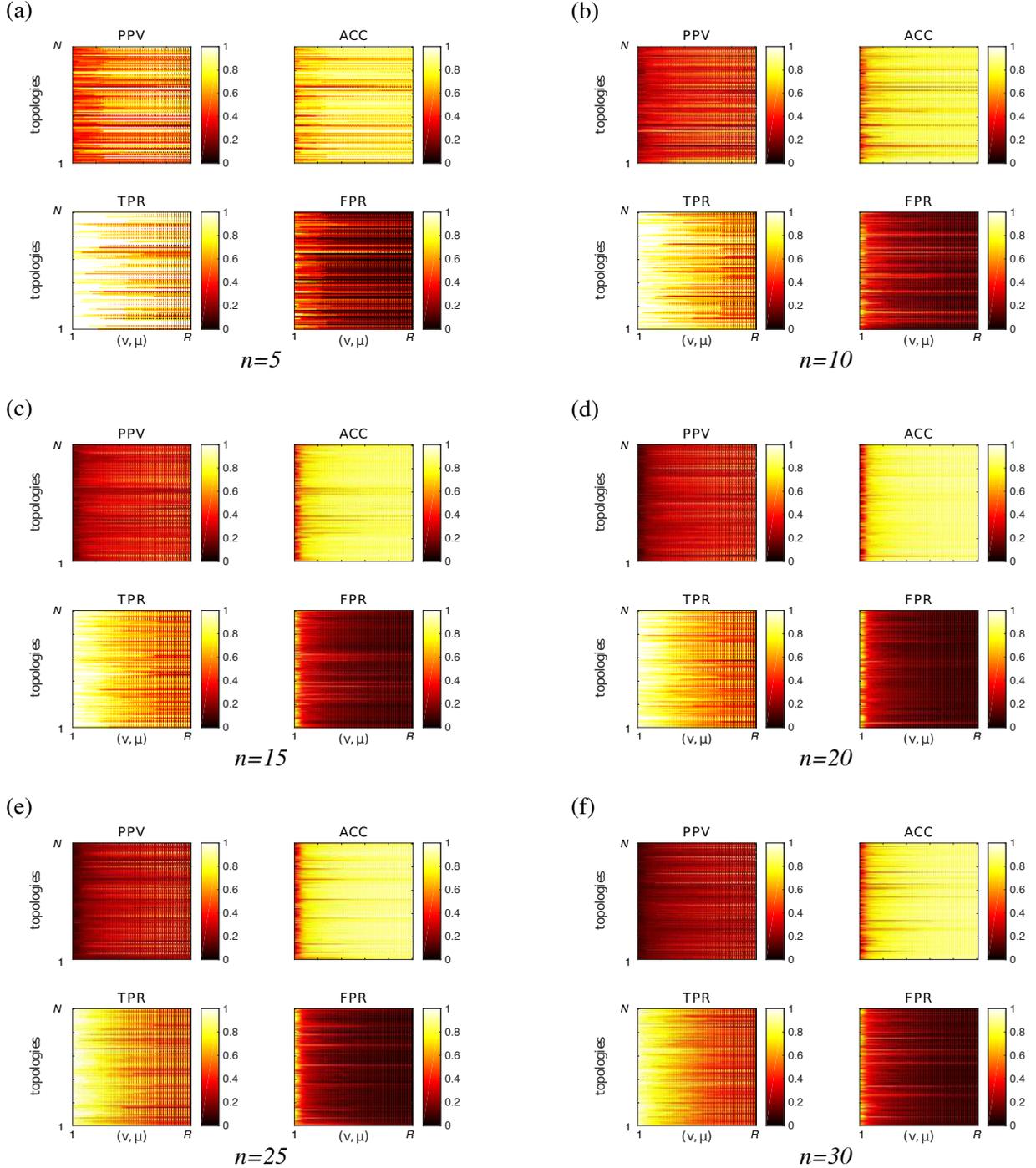}
 \caption{Output of \emph{Step 2}. $4RN$ metrics, with $R=5050$ and $N=100$, are depicted as heat maps, where all the threshold couples $(\nu,\mu) \in [0,0.99]^2$ are represented on the $x$ axis, the topology index is represented on the $y$ axis, and different colors quantify the value of the metrics, with dark (light) colors being representative of low (high) values. The results are shown for known test-topologies with a different number of nodes $n$. (a) $n=5$, (b) $n=10$, (c) $n=15$, (d) $n=20$, (e) $n=25$, (f) $n=30$.}
 \label{fig:S2}
\end{figure}

The output of \emph{Step 2} consists in $4RN$ metrics which can be depicted as heat maps, with different colors quantifying the value of the metrics. Specifically, dark (light) colors are representative of low (high) values (Figure \ref{fig:S2}). As the number of nodes $n$ increases, PPV and TPR take lower values, whilst ACC and FPR do not vary as much.

%
%

\clearpage
\subsubsection*{Output of \emph{Step 3}: acceptable values for thresholds $\nu$ and $\mu$}

\begin{figure}[!h]
 \centering
\includegraphics[scale=.75]{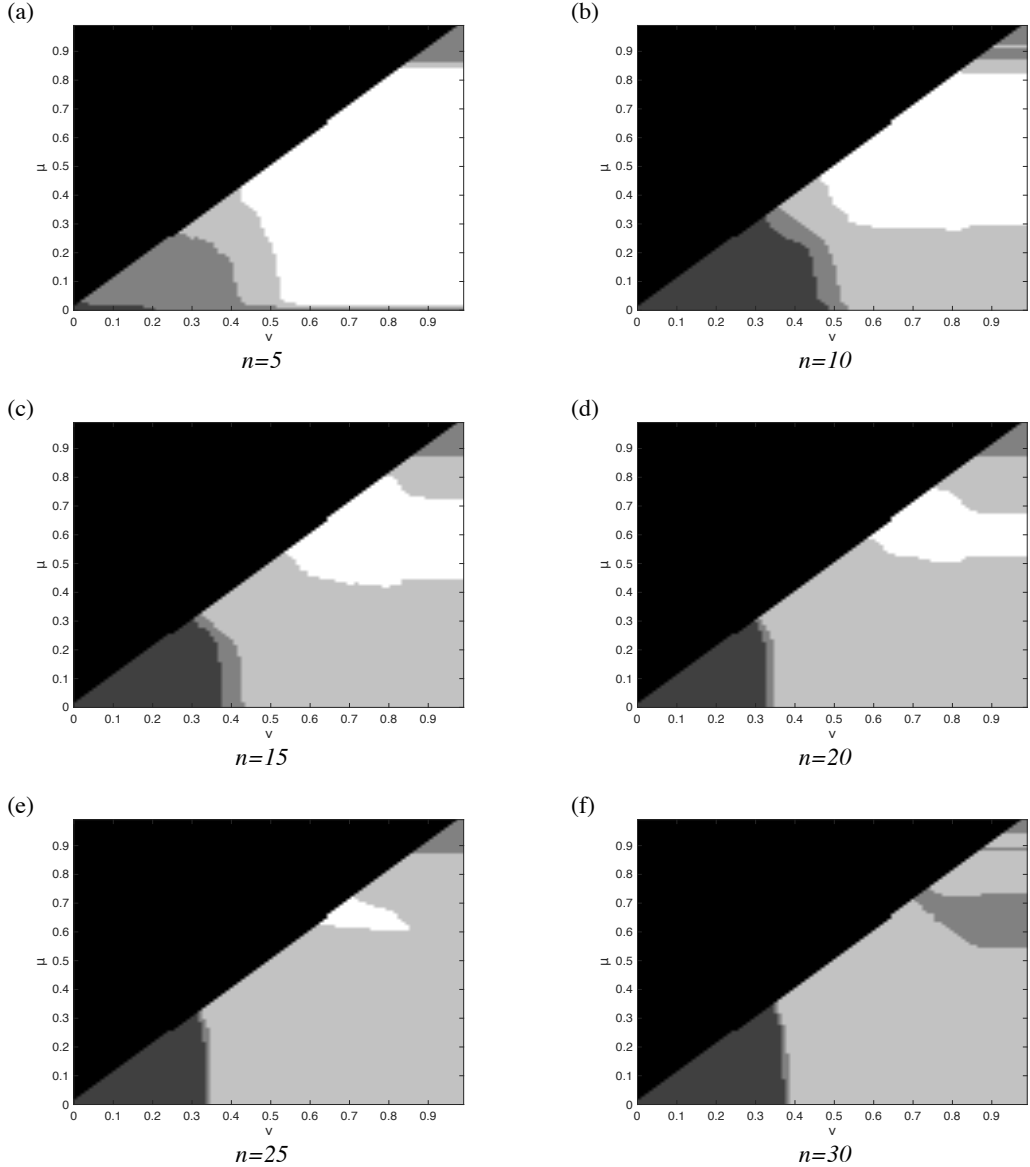}
 \caption{Output of \emph{Step 3}. Acceptable ranges within which $\nu \in [0,0.99]$ and $\mu \in [0,0.99]$ should take values when reconstructing an unknown topology are depicted as a shades-of-gray map, where $\nu$ is represented on the $x$ axis, $\mu$ on the $y$ axis, and with lighter (darker) colors being representative of values for the threshold couples in correspondence to which more (less) metrics bounds conditions in Equation \eqref{eqn:metricsBounds} are verified. White regions correspond to values of $(\nu,\mu)$ for which the bounds conditions are all simultaneously satisfied. The results are shown for known test-topologies with a different number of nodes $n$. (a) $n=5$, (b) $n=10$, (c) $n=15$, (d) $n=20$, (e) $n=25$, (f) $n=30$.}
 \label{fig:S3}
\end{figure}

The output of \emph{Step 3} consists in acceptable ranges within which $\nu$ and $\mu$ should take values when reconstructing an unknown topology. It can be depicted as a shades-of-gray map (Figure \ref{fig:S3}), with lighter (darker) colors being representative of values for the threshold couples in correspondence to which more (less) metrics bounds conditions in Equation \eqref{eqn:metricsBounds} are verified. Ideally, $\nu$ and $\mu$ should be chosen such that they individuate a point belonging to a white region (the bounds conditions are all simultaneously satisfied). If white regions are not found, a sensible choice would be to select threshold values that individuate points belonging to the lightest available region, or to relax the bounds conditions in Equation \eqref{eqn:metricsBounds} and repeat the algorithm.

Note how, as the number of nodes $n$ increases, the area of admissible regions where the threshold conditions in Equation \eqref{eqn:metricsBounds}  are all simultaneously satisfied shrinks (Figure \ref{fig:S3}).

\clearpage
\subsection*{Values of the metrics for different numbers of experiments}
\begin{figure}[!h]
 \centering
\includegraphics[scale=.8]{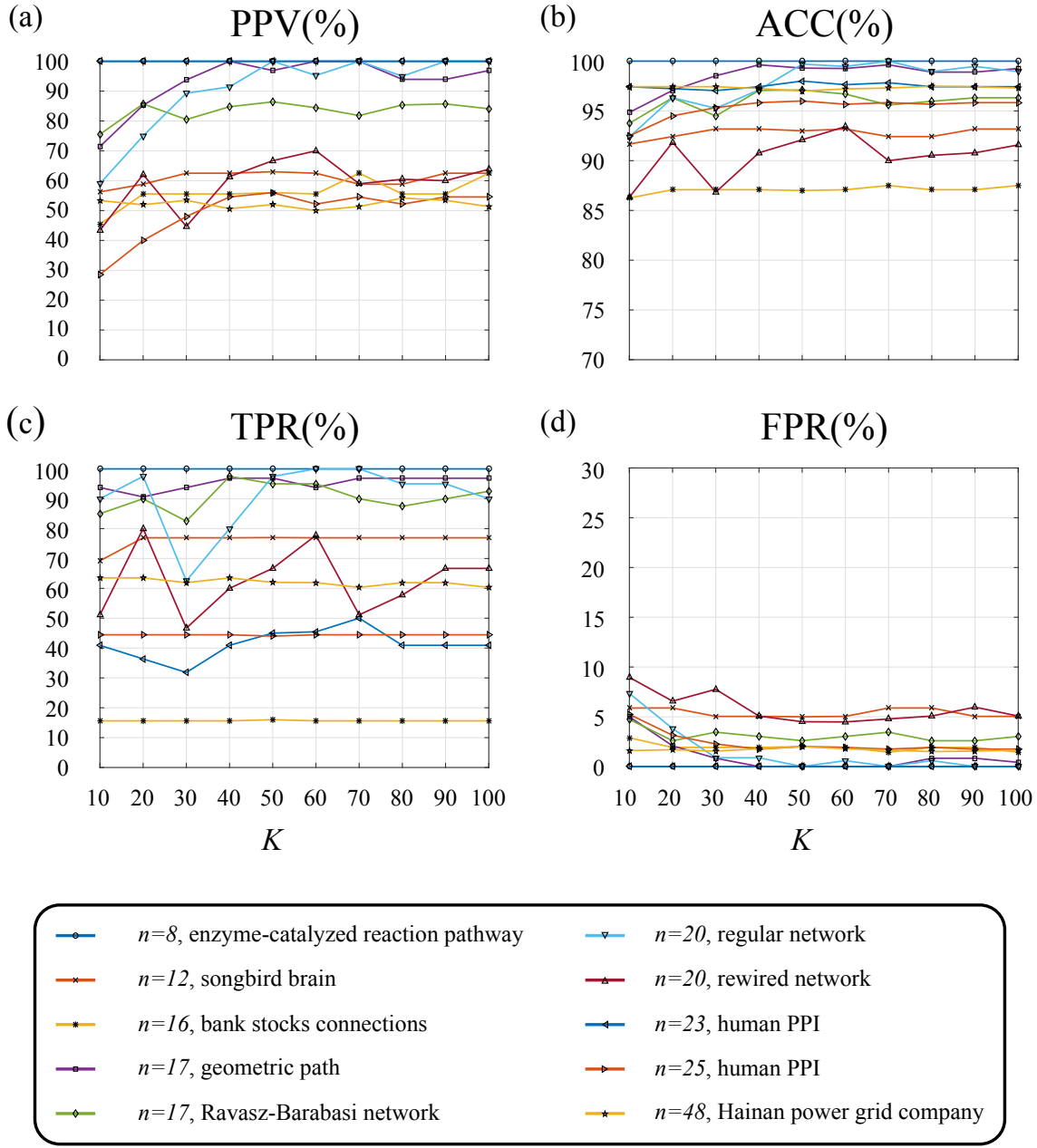}
 \caption{Values of the four metrics for different numbers of experiments $K$. Different colors and symbols refer to different topologies considered in the main text. (a) PPV: Positive Predictive Value, (b) ACC: Accuracy Rate, (c) TPR: True Positive Rate, (d) FPR: False Positive Rate.}
 \label{fig:S4}
\end{figure}

The four standard metrics are quantified in the main text for $K=50$. In Figure \ref{fig:S4} we show their values as a function of the number of experiments $K$ employed to numerically generate the data-set. It is possible to appreciate how, for all the considered topologies, possible fluctuations of the metrics, due to the fact that each experiment is run with random values for both initial conditions and parameters of the nodes, are negligible when $K\ge50$. 

For the sake of simplicity we do not show the value of the metrics for the \emph{4-node} topologies considered in the main text: indeed, for these topologies the values of the metrics do not vary as the number of experiments $K$ increases ($PPV=ACC=TPR=100 \%$, $FPR=0 \%$). 

\end{document}